\documentclass[11pt,a4paper]{article}
\usepackage[top=3cm, bottom=3cm, left=3cm, right=3cm]{geometry}
\usepackage[latin1]{inputenc}
\usepackage{csquotes}
\usepackage{amsmath}
\usepackage{amsthm}
\usepackage{amsfonts}
\usepackage{amssymb}
\usepackage{graphics}
\usepackage{float}
\usepackage[hang,flushmargin]{footmisc} 

\usepackage{amsmath,amsthm,amssymb,graphicx, multicol, array}
\usepackage{enumerate}
\usepackage{enumitem}
\usepackage{xcolor}

\usepackage[pagebackref,bookmarks,colorlinks,breaklinks,linktoc=page]{hyperref}
\hypersetup{linkcolor=blue,citecolor=red,filecolor=blue,urlcolor=blue} 

\usepackage{tabto}

\numberwithin{equation}{section}

\usepackage{tikz}
\usepackage{pgfplots}

\usetikzlibrary{matrix}

\usepackage{mathrsfs}

\DeclareMathOperator{\li}{li}
\DeclareMathOperator{\ord}{ord}

\DeclareMathOperator{\tmod}{mod}

\newtheorem{thm}{Theorem}[section]
\newtheorem{lem}{Lemma}[section]
\newtheorem{theorem}{Theorem}[section]
\newtheorem{lemma}{Lemma}[section]

\newtheorem{dfn}{Definition}[section]

\newcommand{\N}{\mathbb{N}}
\newcommand{\Z}{\mathbb{Z}}
\newcommand{\Q}{\mathbb{Q}}

\newcommand{\C}{\mathbb{C}}
\newcommand{\F}{\mathbb{F}}
\newcommand{\tP}{\mathbb{P}}

\usepackage{fancyhdr}
\usepackage{titletoc}

\titlecontents{section}
[0pt]                                               
{}%
{\contentsmargin{0pt}                               
	\thecontentslabel\enspace%
	\normalsize}
{\contentsmargin{0pt}\normalsize}                        
{\titlerule*[.5pc]{.}\contentspage}                 
[] 
\title{Simple Proof of the Primitive Root Conjecture}
\date{}
\author{N. A. Carella}

\begin{document}
\thispagestyle{empty}
\maketitle

\vskip .25 in
\begin{abstract}
Let \(u\neq \pm 1,v^2\) be a fixed integer, let \(p\geq 2\) be a prime, and let \(\ord_p(u) \mid p-1\) be the multiplicative order of \(u \tmod p\). Define a prime counting function by $\pi(u,x)=\# \{ p\leq x:\ord_p(u)=p-1 \}$. In 1967 Hooley proved a conditional asymptotic formula $\pi(u,x)=\delta(u)x(\log x)^{-1}+O(\log\log x(\log  x)^{-2}$ for the primitive root conjecture. This note proves an unconditional asymptotic formula $\pi(u,x)=\delta(u)x(\log x)^{-1}+O(x(\log  x)^{-2}$ of the same result, where $\delta(u)>0$ is the density constant.  \let\thefootnote\relax\footnote{\today\date{} 
	}
\end{abstract}

\tableofcontents


\section{Introduction} \label{S7776-I}\hypertarget{S7776-I}
The symbols $\N=\{0,1,2,3,\ldots\}$, $\mathbb{Z}=\{0, \pm 1, \pm 2, \ldots \}$  and $\mathbb{P}=\{ 2, 3, 5, \ldots \}$ denote the set of natural numbers, the set integers and the set of prime numbers respectively. Let \(x\geq 1\) be a large real number and let $\pi(x)=\#\{\text{prime }p\leq x\}$ be the prime counting function. The multiplicative order of an element in a finite field $u\in \F_p$ is denoted by $\ord_p u$, see \hyperlink{dfn9955FF.200}{Definition} \ref{dfn9955FF.200}. The constant \(\delta(u)\geq 0\) is the density 
\begin{equation}\label{eq7776.101d}
\delta(u)=\lim_{x\to\infty}\frac{\#\left\{ p\in \mathbb{P}:\ord_p(u)=p-1 \right\}}{\pi(x)}
\end{equation} 
of the subset of primes $\mathscr{P}_u=\left\{ p\in \mathbb{P}:\ord_p(u)=p-1 \right\}$ with a fixed primitive root \(u\) in the set of primes $\tP$.  The expected number of primes \(p\leq x\) with a fixed primitive root \(u\tmod p\) has the asymptotic formula
\begin{equation} \label{eq7776.101f}
\pi _u(x)=\#\left\{ p\leqslant x:\ord_p(u)=p-1 \right\}=\delta(u)\li(x)+O\left(x (\log  x)^{-2} \right),
\end{equation}
where \(\li(x)\) is the logarithm integral, as \(x\rightarrow \infty\). \\

A conditional proof of this result was achieved in \cite{HC1967}, and simplified sketches of the proof appear in {\color{red}\cite[p. 8]{MP2004}}, and similar references.
The determination of the constant \(\delta(u)\geq 0\) for a fixed integer \(u\in \mathbb{Z}\) is an interesting technical problem, {\color{red}\cite[p. 218]{HC1967}}, \cite{LM2014}, \cite{LH1977}, et alii. An introduction to its historical development, and its calculations is covered in {\color{red}\cite[pp. 3--10]{MP2004}}, and \cite{SP2003}. More recent results are given in \cite{SP2021}, \cite{PT2022}, et alia.\\
   
The Artin primitive root conjecture on average
\begin{equation} \label{eq7776.101e}
\frac{1}{x}\sum _{u\leq x} \pi _u(x)=a_0\li(x)+O\left(\frac{x}{(\log  x)^{-c}}\right).
\end{equation} 
The Artin constant is defined by
\begin{eqnarray}\label{eq7776.101k}
a_0=\sum_{n\geq1}\frac{\mu(n}{[\mathcal{K}_n:\Q]}=\prod _{p\geq 2} \left(1-\frac{1}{p(p-1)}\right)
\end{eqnarray}
where $[\mathcal{K}_n:\Q]$ is a number field extension of $\Q$ and \(c>1\) is an arbitrary number, was proved in \cite{SP2003}. These works had shown that almost all admissible integers \(u\) are primitive roots for infinitely many primes. The number of exceptions is a subset of zero density in \(\mathbb{Z}\). The individual quantity \(\pi _u(x)\) in \eqref{eq7776.101f} can be slightly different from the average quantity in \eqref{eq7776.101k}. The variations, discovered by the Lehmers using
numerical experiments, depend on the primes decomposition of the fixed value \(u\), see {\color{red}\cite[p. 220]{HC1967}}, {\color{red}\cite[p. 3]{MP2004}}, and similar references for the exact formula for the density \(\delta(u) \geq 0\). \\

These results have been extended to quadratic numbers fields in \cite{CJ2007}, \cite{NW2000}, \cite{RH2002}, and other number fields in \cite{HJ1986}, \cite{AC2014}, \cite{MS2024}, et cetera. Other related results are given in \cite{MP2007}, \cite{FX2011}, \cite{PS2013}, \cite{AC2014A}, et alii. The Artin primitive root conjecture for functions fields was proved by Bilharz, see {\color{red}\cite[Theorem 10.11]{RN1999}}, and the same conjecture for polynomials over finite fields was proved in \cite{PS1995}. The literature on this topic is vast and extensive.\\

This note introduces a simple unconditional proof of the primitive root conjecture over the integers $\Z$.

\begin{thm} \label{thm7776.100}\hypertarget{thm7776.100} A fixed integer \(u\neq \pm 1,v^2\) is a primitive root mod \(p\) for infinitely many primes \(p\geq 2\). In addition, the density of these primes satisfies 
\begin{equation} \label{eq7776.100i}
\pi _u(x)=\#\left\{ p\leq x:\ord_p(u)=p-1 \right\}=\delta(u) \li(x)+O\left(x (\log  x)^{-2} \right),
\end{equation}
where \(\li(x)\) is the logarithm integral, and \(\delta(u)>0\) is a constant, for all large numbers \(x\geq 1\) unconditionally.
\end{thm} 

\hyperlink{S9955FF}{Section} \ref{S9955FF} to \hyperlink{S6-S}{Section} \ref{S6-S} introduce the notation and some standard results related to or applicable to the investigation of primitive roots in cyclic groups. \hyperlink{S7777-K}{Section} \ref{S7777-K} presents a proof of \hyperlink{thm7776.100}{Theorem} \ref{thm7776.100}.

\section{Characteristic Functions in Finite Fields}\label{S9955FF}\hypertarget{S9955FF}
The characteristic function \(\Psi :G\longrightarrow \{ 0, 1 \}\) of primitive elements is one of the standard analytic tools employed to investigate the various properties of primitive roots in cyclic groups \(G\). Many equivalent representations of the characteristic function $\Psi $ of primitive elements
are possible, a few are investigated here.
\begin{dfn}\label{dfn9955FF.200}\hypertarget{dfn9955FF.200} {\normalfont The class function $ \ord: \F_p ^{\times}\longrightarrow \{d<p:d\mid p-1\}$ defined by $\ord_p u=\min \#\{n<p:u^n\equiv 1 \bmod p\}$ measures the multiplicative order of an element $u\ne0$ in the multiplicative group $\F_p ^{\times}$. An element of maximal multiplicative order $\ord_p u=p-1$ is called a \textit{primitive element}. 
	}
\end{dfn}

\subsection{Divisor Dependent Characteristic Function}		\label{S9955FF-A}\hypertarget{S9955FF-A}
The \textit{divisor-dependent} representation of the characteristic function\index{Divisor dependent characteristic function!primitive root in finite field} detects the multiplicative order of an element by means of the divisors of the totient function $\varphi(p)=p-1$. The version for primitive elements in the prime finite fields $\F_p$ is described below.

\begin{lem} \label{lem9955FF.200DD} \hypertarget{lem9955FF.200DD} Let \(p\geq 2\) be an integer and let \(\chi\) be a multiplicative character of order $\ord  \chi =d$. If \(u\in
	\F_p\) is a nonzero element, then
	\begin{equation}\label{eq9955FF.200DDe}
		\Psi (u)=\frac{\varphi(p-1)}{p-1}\sum _{d\mid \varphi(p)} \frac{\mu(d)}{\varphi(q)}\sum _{\ord \chi =d} \chi(u)
		=\left \{
		\begin{array}{ll}
			1 & \text{   \normalfont if } \ord_p (u)=p-1,  \\
			0 & \text{   \normalfont if } \ord_p (u)\neq p-1, \\
		\end{array} \right .\nonumber
	\end{equation}
	where $\mu:\N\longrightarrow \{-1,0,1\}$ is the Mobius function.
\end{lem}	

There are a few other variant proofs of this result, these are widely available in the literature. Almost every result in the theory of primitive roots in finite fields is based on this characteristic function, but sometimes written in different forms. This technique was developed by Landau, Vinogradov, Davenport, Erdos, and other authors, see \cite{LE1927}, \cite{VI1930}, \cite{DH1937}, \cite{EP1945}, \cite{WR2001}, et alii. An extension of this characteristic function to the finite ring $\Z/n\Z$ is presented in {\color{red}\cite[Lemma 4]{MG1998}}. The main obstacle in this technique is the primes decomposition of the totients $p-1$. Specifically, the finite sum $\sum _{d\mid p-1} |\mu(d)|\ll p^{\varepsilon}$ is an obstacle in many applications.
\subsection{Divisorfree Characteristic Functions}\label{S9955FF-B} \hypertarget{S9955FF-B}
A new \textit{divisorfree} representation of the characteristic function\index{Divisorfree characteristic function!primitive root in finite field} of primitive element is introduced in this section. It detects the order \(\text{ord}_p
(u) \geq 1\) of the element \(u\in \mathbb{F}_p\) by means of the solutions of the equation \(\tau ^s-u=0\) in \(\mathbb{F}_p\), where
\(u,\tau\) are constants, and $s\in \mathscr{R}=\{s<p:\gcd(s,p-1)=1\}$ is a variable. 

\begin{lem} \label{lem9955FF.200DF} \hypertarget{lem9955FF.200DF}  Let \(p\geq 2\) be a prime and let \(\tau\) be a primitive root mod \(p\) and  let \(\psi \neq 1\) be a nonprincipal additive character of order $\ord  \psi =p$. If \(u\in
	\mathbb{F}_p\) is a nonzero element, then
	\begin{equation}
		\Psi (u)=\sum _{\substack{1\leq s\leq p-1\\\gcd (s,p-1)=1}} \frac{1}{p}\sum _{0\leq t\leq p-1} \psi \left ((\tau ^s-u)t\right)
		=\left \{
		\begin{array}{ll}
			1 & \text{   \normalfont if } \ord_p (u)=p-1,  \\
			0 & \text{   \normalfont if } \ord_p (u)\neq p-1. \\
		\end{array} \right .\nonumber
	\end{equation}
\end{lem}	

\begin{proof}[\textbf{Proof}] Set the additive character $\psi(s) =e^{i 2\pi  as/p}\in \C$. Now, as the index $s\in \mathscr{R}$ ranges over the integers relatively prime to $\varphi(p)=p-1$, the element $\tau ^n\in \F_p ^{\times}$ ranges over the primitive roots modulo $p$. Accordingly, the equation 
\begin{eqnarray}\label{eq9977FF.200DF2}
a=\tau ^n- u=0
\end{eqnarray}has a unique solution $s\in \mathscr{R}$ if and only if the fixed element $u\in \F_p$ is a primitive root. This implies that the inner sum in 	
	\begin{equation}\label{eq9977FF.200DF}
		\sum _{\substack{1\leq s\leq p-1\\\gcd (s,p-1)=1}} \frac{1}{p}\sum _{0\leq t< p} e^{i 2\pi \frac{(\tau ^s-u)t}{p}}=
		\left \{\begin{array}{ll}
			1 & \text{   \normalfont if } \ord_{p} (u)=p-1,  \\
			0 & \text{   \normalfont if } \ord_{p} (u)\ne p-1. \\
		\end{array} \right.
	\end{equation} 
	collapses to $\sum _{0\leq s< p} e^{i 2\pi as/p}=\sum _{0\leq s< p} 1=p $. Otherwise, if the element $u\in \F_p$ is not a primitive root, then the equation \eqref{eq9977FF.200DF2} has no solution $n\in \mathscr{R}$, and the inner sum in \eqref{eq9977FF.200DF} collapses to $\sum _{0\leq s< p} e^{i 2\pi as/p}=0$,
	this follows from the geometric series formula $\sum_{0\leq n\leq  N-1} w^n =(w^N-1)/(w-1)$, where $w=e^{i 2\pi a/(p-1)}\ne1$ and $N=p$. 
	This completes the verification.	 
\end{proof}

Another approach to the construction of divisorfree characteristic function is via the discrete logarithm $\log_{\tau}:\F_p^{\times}\longrightarrow \F_p$ in finite field $\F_p$.

\begin{lem} \label{lem9955FF.200DFL} \hypertarget{lem9955FF.200DFL}  Let \(p\geq 2\) be a prime and let \(\tau\) be a primitive root mod \(p\) and  let \(\psi \neq 1\) be a nonprincipal additive character of order $\ord  \psi =p$. If \(u\in
	\mathbb{F}_p\) is a nonzero element, then
	\begin{equation}
		\Psi (u)=\sum _{\substack{1\leq s< p-1\\\gcd (s,p-1)=1}} \frac{1}{p}\sum _{0\leq t< p} \psi \left ((s-\log_{\tau}u)t\right)
		=\left \{
		\begin{array}{ll}
			1 & \text{   \normalfont if } \ord_p (u)=p-1,  \\
			0 & \text{   \normalfont if } \ord_p (u)\neq p-1. \\
		\end{array} \right .\nonumber
	\end{equation}
\end{lem}	

\begin{proof}[\textbf{Proof}] Let $\mathscr{R}=\{s<p:\gcd(s,p-1)=1\}$ and set the additive character $\psi(s) =e^{i 2\pi  as/p}\in \C$. Consider a fixed element $u=\tau^s\in \F_p$ and its discrete logarithm $\log_{\tau}u=s$. As the index varies over the set $s\in \mathscr{R}$ ranges over the integers relatively prime to $\varphi(p)=p-1$, the equation 
	\begin{equation}\label{eq9977FF.300DFL1}
		a=s- \log_{\tau}u=0
	\end{equation} 
	has a unique solution $s=\log_{\tau}u$ if and only if the fixed element $u\in \F_p$ is a primitive root. This implies that the inner sum in 	
	\begin{equation}\label{eq9977FF.300DFL2}
		\sum _{\substack{1\leq s< p-1\\\gcd (s,p-1)=1}} \frac{1}{p}\sum _{0\leq t< p} e^{i 2\pi \frac{(s-\log_{\tau}u)t}{p}}=
		\left \{\begin{array}{ll}
			1 & \text{   \normalfont if } \ord_{p} (u)=p-1,  \\
			0 & \text{   \normalfont if } \ord_{p} (u)\ne p-1. \\
		\end{array} \right.
	\end{equation} 
	collapses to $\sum _{0\leq s< p} e^{i 2\pi as/p}=\sum _{0\leq s< p-1} 1=p $. Otherwise, if the element $u-1\in \F_p$ is not a primitive root, then the equation \eqref{eq9977FF.300DFL1} has no solution $s\in \mathscr{R}$ , and the inner sum in \eqref{eq9977FF.300DFL2} collapses to $\sum _{0\leq s< p-1} e^{i 2\pi as/p}=0$,
	this follows from the geometric series formula $\sum_{0\leq n\leq  N-1} w^n =(w^N-1)/(w-1)$, where $w=e^{i 2\pi a/p}\ne1$ and $N=p$. 
	This completes the verification.	 
\end{proof}
\section{Finite Summation Kernel and Gauss Sum}
An upper bound for some elementary exponential sums are provided in this section. 
\subsection{Finite Summation Kernel}
\begin{lemma}   \label{lem5555.400B}\hypertarget{lem5555.400B}  Let \(p\geq 2\) be large prime, and let $\omega=e^{i2 \pi/p} $ be a $p$th root of unity. Then,
\begin{equation}
		\sum_{1 \leq t\leq p-1}\Bigg | \sum_{\substack{1\leq n\leq p-1\\\gcd(n,p-1)=1}} \omega^{tn}  \Bigg  |\ll  p^{1+\delta} \log p ,
\end{equation}		where $\delta>0$ is a small number. 
\end{lemma} 

\begin{proof}[\textbf{Proof}] Use the inclusion exclusion principle to rewrite the exponential sum as
	\begin{eqnarray} \label{eq5555.400i}
		\sum_{\substack{1\leq n\leq p-1\\\gcd(n,p-1)=1}}\omega^{tn}&=& \sum_{n \leq p-1} \omega^{tn}  \sum_{\substack{d \mid p-1 \\ d \mid n}}\mu(d)  \nonumber \\
		&=& \sum_{d \mid p-1} \mu(d) \sum_{\substack{n \leq p-1 \\ d \mid n}} \omega^{tn}\nonumber \\
		& =&\sum_{d\mid p-1} \mu(d) \sum_{m \leq (p-1)/ d} \omega^{dtm} \\
		&=& \sum_{d \mid p-1} \mu(d) \frac{\omega^{dt}-\omega^{dt((p-1)/d+1)}}{1-\omega^{dt}} \nonumber.
	\end{eqnarray} 
	Now, the parameters are $p$ prime, $\omega=e^{i2 \pi/p}$, the integers $t \in [1, p-1]$ and $d \leq p-1<p$. This data implies that $\pi dt/p\ne k \pi $ with $k \in \mathbb{Z}$, so the sine function $\sin(\pi dt/p)\ne 0$ is well defined. Consequently, the absolute value satisfies
	\begin{equation}
		\left |\frac{\omega^{dt}-\omega^{dt((p-1)/d+1)}}{1-\omega^{dt}} \right |\leq 	\left | \frac{2}{\sin( \pi dt/ p)} \right |.
	\end{equation}
For each $d\mid p-1$, the map $t\longrightarrow z\equiv dt \bmod p$ is a permutation in the finite field $\F_p$. Thus, using standard manipulations, and $z/2 \leq \sin(z) <z$ for $0<|z|<\pi/2$, the last expression becomes
	\begin{eqnarray}
		\sum_{1 \leq t\leq p-1}\Bigg | 	\sum_{\substack{1\leq n\leq p-1\\\gcd(n,p-1)=1}}\omega^{tn}\Bigg  |&\leq&\sum_{d \mid p-1,} \sum_{1 \leq t\leq p-1}	\left | \frac{2}{\sin( \pi dt/ p)} \right |\\
		&\leq&\sum_{d \mid p-1,} \sum_{1 \leq z\leq p-1}	 \frac{2p}{ \pi z} \nonumber\\
		&\ll&  p^{1+\delta} \log p \nonumber,
	\end{eqnarray}
	where $\sum_{d \mid p-1}1=d(p-1)\ll p^{\delta}$ is the number of divisor in $p-1$ and $\delta>0$ is a small number. 
\end{proof}

\subsection{Gauss Sum}
Some elementary exponential sums estimates are provided in this section. 
\begin{lemma}   \label{lem1234A.150A}\hypertarget{lem1234A.150A}  
	{\normalfont (Gauss sums)} Let \(p\geq 2\) be a prime, let $\chi(t)=e^{i2 \pi t/p} $ and  $\psi(t)=e^{i2\pi  \tau^t/p}$ be a pair of characters. Then, the Gaussian sum has the upper bound
	\begin{equation} \label{eq3-355}
		\left |\sum_{1 \leq t \leq p-1}    \chi(t) \psi(t) \right | \leq 2 p^{1/2} \log p.\nonumber
	\end{equation}
	
\end{lemma} 

\section{Estimates of Power Exponential Sums} \label{S9933Q-A}\hypertarget{S9933Q-A}
The estimate for the power sum with relatively prime index is based on the identity
\begin{equation}\label{eq9933Q.210c}
\frac{1}{p} \sum_{0 \leq t\leq p-1,}  \sum_{0 \leq s\leq p-1} \omega^{t(n-s)}f(s)=f(n),\end{equation}
where $\omega=e^{i2\pi/p}$ and $x  \leq p -1$. The case of interest here is the exponential sum
	\begin{equation} \label{eq9933Q.346b}
		\sum_{\substack{1\leq n\leq p-1\\\gcd(n,p-1)=1}} e^{\frac{i2\pi b \tau^n}{p}}= \sum_{\substack{1\leq n\leq p-1\\\gcd(n,p-1)=1}}\frac{1}{p} \sum_{0 \leq t\leq p-1,}  \sum_{1 \leq s\leq p-1} \omega^{t(n-s)}e^{\frac{i2\pi b \tau^s}{p}}.
	\end{equation}
Use the inclusion exclusion principle to rewrite the exponential sum as
	\begin{equation}\label{eq9933Q.346d}
		\sum_{\substack{1\leq n\leq p-1\\\gcd(n,p-1)=1}} e^{ \frac{i2\pi b \tau^n}{p}} 
		= \sum_{1\leq  n \leq p-1}\frac{1}{p} \sum_{0 \leq t\leq p-1,}  \sum_{1 \leq s\leq p-1} \omega^{t(n-s)}e^{\frac{i2\pi b \tau^s}{p}} \sum_{\substack{d \mid p-1 \\ d \mid n}}\mu(d)   .
	\end{equation} 
Substituting $t=0$ in \eqref{eq9933Q.210c} contributes $-\varphi(p-1)/p$, and rearranging it yield
	\begin{eqnarray}\label{eq9933Q.346f}
		&& \sum_{\substack{1\leq n\leq p-1\\\gcd(n,p-1)=1}} e^{ \frac{i2\pi b \tau^n}{p}} \\
		&=& \sum_{1\leq  n \leq p-1}\frac{1}{p} \sum_{1 \leq t\leq p-1,}  \sum_{1 \leq s\leq p-1} \omega^{t(n-s)}e^{\frac{i2\pi b \tau^s}{p}} \sum_{\substack{d \mid p-1 \\ d \mid n}}\mu(d) -\frac{\varphi(p-1)}{p} \nonumber \\
		&=&\frac{1}{p} \sum_{1 \leq t\leq p-1} \left ( \sum_{1 \leq s\leq p-1} \omega^{-ts}e^{\frac{i2\pi b \tau^s}{p}}\right )\left (\sum_{d \mid p-1} \mu(d) \sum_{\substack{1\leq n \leq p-1, \\ d \mid n}}   \omega^{tn} \right ) -\frac{\varphi(p-1)}{p} \nonumber.
	\end{eqnarray}

\begin{theorem}  \label{thm9933Q.346}\hypertarget{thm9933Q.346}  Let \(p\geq 2\) be a large prime, and let $\tau $ be a primitive root modulo $p$. Then,
	\begin{equation}
		\sum_{\substack{1\leq n\leq p-1\\\gcd(n,p-1)=1}} e^{i2\pi b \tau^n/p} \ll  p^{1/2+\delta}(\log p)^2 \nonumber,
		\end{equation} 
		where $\delta>0$ is a small real number and the implied constant is independent of $b\ne0$. 	
\end{theorem}

\begin{proof}[\textbf{Proof}] There are many proofs of this result, consult \cite{ML1972}, {\color{red}\cite[Theorem 6]{FS2000}}, and related results are given in \cite{FS2001}, \cite{GM2005}, \cite{CC2009}, and {\color{red}\cite[Theorem 1]{GK2005}}. The upper bound given in \hyperlink{thm9933Q.346}{Theorem} \ref{thm9933Q.346} seems to be optimum. A different proof, which has a weaker upper bound, appears in {\color{red}\cite[Theorem 6]{FS2000}}, and related results are given in \cite{CC2009}, \cite{FS2002}, \cite{GK2005}, and {\color{red}\cite[Theorem 1]{GK2005}}.\end{proof}

\subsection{Finite Fourier Transform of the Power Exponential Sum} 
For any fixed $ 0 \ne b \in \mathbb{F}_p$, the map $ \tau^n \longrightarrow b \tau^n$ is one-to-one (permutation) in $\mathbb{F}_p$. Consequently, the subsets 
\begin{equation} \label{eq9933RPI.500b}
	\{ \tau^n: \gcd(n,p-1)=1 \}\quad \text { and } \quad  \{ b\tau^n: \gcd(n,p-1)=1 \} \subset \mathbb{F}_p
\end{equation} have the same cardinalities. As a direct consequence the exponential sums 
\begin{equation} \label{eq9933RPI.500d}
	V(b)=\sum_{\substack{1\leq n\leq p-1\\\gcd(n,p-1)=1}}e^{i2\pi b \tau^n/p} \quad \text{ and } \quad V(1)=\sum_{\substack{1\leq n\leq p-1\\\gcd(n,p-1)=1}} e^{i2\pi \tau^n/p},
\end{equation}
have the same upper bound up to an error term. An asymptotic relation for the finite Fourier transform (FFT) of the exponential sums \eqref{eq9933RPI.500d} is provided here. 

\begin{theorem}   \label{thm9933ERP.220V}\hypertarget{thm9933ERP.220V}  Let \(p\geq 2\) be a large prime. If $\tau $ be a primitive root modulo $p$ and $a<x=o(p)$ is not a primitive root, then
	\begin{equation} 
	\widehat{V(a)}=	\sum_{1\leq b\leq  p-1}	 e^{-i2\pi \frac{ab}{p}}	\sum_{\substack{1\leq n\leq p-1\\\gcd(n,p-1)=1}} e^{\frac{i2\pi b \tau^n}{p}} =-  \sum_{\substack{1\leq n\leq p-1\\\gcd(n,p-1)=1}} e^{\frac{i2\pi  \tau^n}{p}} + O(p^{1/2+\delta} (\log p)^2)\nonumber,
	\end{equation} 
	where $\delta>0$ is a small number and the implied constant is independent of $ a,b \in [1, p-1]$. 	
\end{theorem} 
\begin{proof}[\textbf{Proof}] For $a\in[1,x]$ and $b\in[1,p-1]$, the exponential sum has the representation 
	\begin{eqnarray} \label{eq9933RPI.500f}
		V(b)&=& \sum_{\substack{1\leq n\leq p-1\\\gcd(n,p-1)=1}} e^{\frac{i2\pi b \tau^n}{p}} \\
		&=&\frac{1}{p} \sum_{1 \leq t\leq p-1} \left ( \sum_{1 \leq s\leq p-1} \omega^{-ts}e^{\frac{i2\pi b \tau^s}{p}}\right )\left (\sum_{d \mid p-1} \mu(d) \sum_{\substack{n \leq p-1, \\ d \mid n}}   \omega^{tn} \right ) -\frac{\varphi(p-1)}{p}\nonumber,
	\end{eqnarray} 
	confer equations \eqref{eq9933Q.346b}, \eqref{eq9933Q.346d} and \eqref{eq9933Q.346f} for more details. In particular, for $b=1$, 
	\begin{eqnarray} \label{eq9933RPI.500h}
		V(1)&=& 	\sum_{\substack{1\leq n\leq p-1\\\gcd(n,p-1)=1}}e^{\frac{i2\pi  \tau^n}{p}} \\
		&=& \frac{1}{p} \sum_{1 \leq t\leq p-1} \left ( \sum_{1 \leq s\leq p-1} \omega^{-ts}e^{\frac{i2\pi a \tau^s}{p}}\right )\left (\sum_{d \mid p-1} \mu(d) \sum_{\substack{n \leq p-1, \\ d \mid n}}   \omega^{tn} \right ) -\frac{\varphi(p-1)}{p}\nonumber,
	\end{eqnarray}
	respectively. Differencing \eqref{eq9933RPI.500f} and \eqref{eq9933RPI.500h} produces 
	\begin{eqnarray} \label{eq9933RPI.500i}
		V(b)-V(1)&= &	\sum_{\substack{1\leq n\leq p-1\\\gcd(n,p-1)=1}} e^{\frac{i2\pi b \tau^n}{p}} -\sum_{\substack{1\leq n\leq p-1\\\gcd(n,p-1)=1}} e^{\frac{i2\pi  \tau^n}{p}} \\
		&=&     \frac{1}{p} \sum_{1 \leq t\leq p-1} \left ( \sum_{1 \leq s\leq p-1} \omega^{-ts}e^{\frac{i2\pi  b \tau^s}{p}}-\sum_{1 \leq s\leq p-1} \omega^{-ts}e^{\frac{i2\pi  \tau^s}{p}}\right ) \nonumber \\
		&& \times \left (\sum_{d \mid p-1} \mu(d) \sum_{\substack{n \leq p-1, \\ d \mid n}}   \omega^{tn} \right ) \nonumber.
	\end{eqnarray}
	Taking the finite Fourier transform of the difference $D(b)=V(b)-V(1)$ returns 

	\begin{eqnarray} \label{eq9933RPI.500j}
		\widehat{D(a)}&=&	\sum_{1\leq b\leq  p-1}	 e^{-i2\pi \frac{ab}{p}}\left( \sum_{\substack{1\leq n\leq  p-1\\\gcd(n,p-1)=1}} e^{\frac{i2\pi b \tau^n}{p}} -\sum_{\substack{1\leq n\leq  p-1\\\gcd(n,p-1)=1}} e^{\frac{i2\pi  \tau^n}{p}}\right)  \\
		&=&  \frac{1}{p} \sum_{1\leq b\leq  p-1}	 e^{-i2\pi \frac{b}{p}}  \sum_{1 \leq t\leq p-1} \left ( \sum_{1 \leq s\leq p-1} \omega^{-ts}e^{\frac{i2\pi b \tau^s}{p}}-\sum_{1 \leq s\leq p-1} \omega^{-ts}e^{\frac{i2\pi a \tau^s}{p}}\right ) \nonumber \\
		&&\hskip 1.75in \times \left (\sum_{d \mid p-1} \mu(d) \sum_{\substack{n \leq p-1, \\ d \mid n}}   \omega^{tn} \right ) \nonumber\\
		&=&   \frac{1}{p} \sum_{1 \leq t\leq p-1} \left ( \sum_{1 \leq s\leq p-1} \omega^{-ts}\sum_{1\leq b\leq  p-1}	  e^{\frac{i2\pi b (\tau^s-a)}{p}}\right.  \nonumber \\
		&&\hskip .15in-\left .\sum_{1\leq b\leq  p-1}	 e^{-i2\pi \frac{ab}{p}}   \sum_{1 \leq s\leq p-1} \omega^{-ts}e^{\frac{i2\pi  \tau^s}{p}}\right ) \times \left (\sum_{d \mid p-1} \mu(d) \sum_{\substack{n \leq p-1, \\ d \mid n}}   \omega^{tn} \right ) \nonumber.
	\end{eqnarray}
	Now in the range $a<x=o(p)$, $\tau^s-a\ne0$ for any $s\in\{s<p-1:\gcd(s,p-1)=1\}$. Thus, using the geometric sum identity $\sum_{1\leq u\leq  p-1}	 e^{i2\pi au/p}=-1$ to simplify the last expression yields
	\begin{eqnarray} \label{eq9933RPI.500l}
		\widehat{D(a)}&=&	\sum_{1\leq b\leq  p-1}	 e^{-i2\pi \frac{ab}{p}}\left( \sum_{\substack{1\leq n\leq  p-1\\\gcd(n,p-1)=1}} e^{\frac{i2\pi b \tau^n}{p}} -\sum_{\substack{1\leq n\leq  p-1\\\gcd(n,p-1)=1}} e^{\frac{i2\pi  \tau^n}{p}}\right)  \nonumber\\
		&=&   \frac{1}{p} \sum_{1 \leq t\leq p-1} \left ( (-1)(-1)-(-1)  \sum_{1 \leq s\leq p-1} \omega^{-ts}e^{\frac{i2\pi  \tau^s}{p}}\right ) \nonumber \\
		&& \times \left (\sum_{d \mid p-1} \mu(d) \sum_{\substack{n \leq p-1, \\ d \mid n}}   \omega^{tn} \right ) .
	\end{eqnarray}
	Rearranging the last equation yield
	\begin{eqnarray} \label{eq9933RPI.500k}
		\widehat{V(a)}&=&\sum_{1\leq b\leq  p-1}	 e^{-i2\pi \frac{ab}{p}}\sum_{\substack{1\leq n\leq  p-1\\\gcd(n,p-1)=1}} e^{\frac{i2\pi b \tau^n}{p}} \\
		&=& -\sum_{\substack{1\leq n\leq  p-1\\\gcd(n,p-1)=1}} e^{\frac{i2\pi  \tau^n}{p}}  +  \frac{1}{p} \sum_{1 \leq t\leq p-1} \left ( 1-  \sum_{1 \leq s\leq p-1} \omega^{-ts}e^{\frac{i2\pi  \tau^s}{p}}\right ) \nonumber \\
		&&\hskip 2.5in \times \left (\sum_{d \mid p-1} \mu(d) \sum_{\substack{n \leq p-1, \\ d \mid n}}   \omega^{tn} \right ) \nonumber\\
		&=&-\sum_{\substack{1\leq n\leq  p-1\\\gcd(n,p-1)=1}} e^{\frac{i2\pi  \tau^n}{p}}+\widehat{R(a)} \nonumber.
	\end{eqnarray}
	
	By \hyperlink{lem5555.400B}{Lemma} \ref{lem5555.400B}, the relatively prime summation kernel is bounded by
	\begin{eqnarray} \label{eq9933RPI.500m}
		\sum_{1 \leq t\leq p-1}	\Bigg |\sum_{d \mid p-1} \mu(d) \sum_{\substack{n \leq p-1, \\ d \mid n}}   \omega^{tn} \Bigg | 
		&=& \sum_{1 \leq t\leq p-1}\Bigg | \sum_{\gcd(n, p-1)=1}\omega^{tn} \Bigg |  \\ 
		&\ll &  p^{1+\delta}\log p\nonumber, 
	\end{eqnarray}
	where $\delta>0$ is a small number and by \hyperlink{lem1234A.150A}{Lemma} \ref{lem1234A.150A}, the difference including Gauss sum is bounded by
	\begin{eqnarray} \label{eq9933RPI.500o}
		\Bigg | 1-  \sum_{1 \leq s\leq p-1} \omega^{-ts}e^{\frac{i2\pi  \tau^s}{p}}\Bigg |=	\Bigg | 1- \sum_{1 \leq s\leq p-1} \chi(s) \psi(s) \Bigg| 
		&\leq & 2 p^{1/2} \log p, 
	\end{eqnarray}
	where  $\chi(s)=e^{i \pi s t/p}$, and $ \psi(s)=e^{i2\pi a \tau^s/p}$. Taking absolute value of the remainder term
	in (\ref{eq9933RPI.500k}) and replacing (\ref{eq9933RPI.500m}), and  (\ref{eq9933RPI.500o}), return
	\begin{eqnarray} \label{eq9933RPI.500p}
		|\widehat{R(a)}|	&=&\frac{1}{p} \Bigg |\sum_{1 \leq t\leq p-1} \Bigg ( 1-  \sum_{1 \leq s\leq p-1} \omega^{-ts}e^{\frac{i2\pi  \tau^s}{p}}\Bigg ) \cdot\Bigg (\sum_{d \mid p-1} \mu(d) \sum_{\substack{n \leq p-1, \\ d \mid n}}   \omega^{tn} \Bigg ) \Bigg | \nonumber\\
		&=&\frac{1}{p}\sum_{1 \leq t\leq p-1}\Bigg | 1-  \sum_{1 \leq s\leq p-1} \omega^{-ts}e^{\frac{i2\pi  \tau^s}{p}} \Bigg | \cdot \Bigg | \sum_{d \mid p-1} \mu(d) \sum_{\substack{n \leq p-1, \\ d \mid n}}   \omega^{tn} \Bigg | \nonumber\\
		&\ll &\frac{1}{p}(2 p^{1/2} \log p)\cdot \sum_{1 \leq t\leq p-1}\Bigg | \sum_{d \mid p-1} \mu(d) \sum_{\substack{n \leq p-1, \\ d \mid n}}   \omega^{tn} \Bigg | \nonumber\\
		&\ll &\frac{1}{p}(2 p^{1/2} \log p)\cdot (p^{1+\delta} \log p) \nonumber\\[.3cm]
		&\ll & p^{1/2+\delta} (\log p)^2,
	\end{eqnarray}
	where the implied constant depends on the number of divisors of $p-1$ and it is independent of $a,b\in [1,p-1]$.
\end{proof}


\begin{lem}   \label{lem9933ERP.230V}\hypertarget{lem9933ERP.230V}  Let \(p\geq 2\) be a large prime. If $\tau $ be a primitive root modulo $p$ and $a<x=o(p)$ is not a primitive root, then
	\begin{equation} 
	\Bigg|\widehat{V(a)}\Bigg|=	\Bigg|\sum_{1\leq b\leq  p-1}	 e^{-i2\pi \frac{ab}{p}}	\sum_{\substack{1\leq n\leq p-1\\\gcd(n,p-1)=1}} e^{\frac{i2\pi b \tau^n}{p}}\Bigg| = O(p^{1/2+\delta} (\log p)^2)\nonumber,
	\end{equation} 
where $\delta>0$ is a small number and the implied constant is independent of $ a,b \in [1, p-1]$. 	
\end{lem} 
\begin{proof}[\textbf{Proof}] The second line in the estimation of the upper bound in \eqref{eq9933ERP.230d} follows from \hyperlink{thm9933ERP.220V}{Theorem} \ref{thm9933ERP.220V} and the fourth line follows from \hyperlink{thm9933Q.346}{Theorem} \ref{thm9933Q.346}: 
\begin{eqnarray} \label{eq9933ERP.230d}
\bigg|\widehat {V(a)}\bigg |&=& 	\Bigg|\sum_{1\leq b\leq  p-1}	 e^{-i2\pi \frac{ab}{p}}	\sum_{\substack{1\leq n\leq p-1\\\gcd(n,p-1)=1}} e^{\frac{i2\pi b \tau^n}{p}}\Bigg| \\[.2cm]
&=&\left |- \sum_{\substack{1\leq n\leq p-1\\\gcd(n,p-1)=1}} e^{i2\pi \frac{\tau ^n}{p}}+ O(p^{1/2+\delta} \log^2 p)  \right | \nonumber\\[.2cm]
&\ll &\left |\sum_{\substack{1\leq n\leq p-1\\\gcd(n,p-1)=1}} e^{i2 \pi \frac{\tau ^n}{p}} \right |+p^{1/2+\delta} (\log p)^2  \nonumber\\[.4cm]
&\ll&  p^{1/2+\delta} (\log p)^2\nonumber,
\end{eqnarray}
where $\delta>0$ is a small number. 
\end{proof}
\section{Evaluation of the Main Term} \label{S7776-MT}\hypertarget{S7776-MT}
The asymptotic formula for the main term occurring in the proof of \hyperlink{thm7776.100}{Theorem} \ref{thm7776.100} is evaluated in this section.\\

\begin{lem} \label{lem7776.200}\hypertarget{lem7776.200} {\normalfont {(\color{red}\cite[Lemma 1]{SP1969}}) } Let \(x\geq 1\) be a large number, and let \(\varphi (n)\) be the Euler totient function. Then
\begin{equation}
\sum_{p\leq x }\frac{\varphi(p-1)}{p-1} =a_1\li(x) +
	O\left(\frac{x}{(\log x)^c}\right) ,
	\end{equation}
where \(\li(x)\) is the logarithm integral, $a_1=0.373955\ldots $, and \(c> 1\) is an arbitrary constant, as \(x \rightarrow \infty\). 
\end{lem}

A more general version of this Lemma is proved in \cite{VR1973}, and related discussions are given in {\color{red}\cite[p.\ 16]{MP2004}}. The generalization to number fields appears in \cite{HJ1984}. These results are ubiquitous in various results in Number Theory. Here, the logarithm integral is defined by
\begin{equation}  \label{eq7776.200c}
\li(x)=\int_2^x \frac{1}{\log t} dt=\frac{x}{\log x} +O\left ( \frac{x}{(\log x)^2} \right).
\end{equation} \\
In this application, the constant  
\begin{equation} \label{eq7776.200f}
a_1=\prod_{p \geq 2 } \left(1-\frac{1}{p(p-1)}\right) 
\end{equation}  
coincides to the average density of primitive roots modulo $p$. The average density of primitive roots was proved in \cite{SP1969}, (the heuristic is due to Artin, and the average density is known as Artin constant). \\

\begin{lem}  \label{lem7776.200B}\hypertarget{lem7776.200B}  Let \(x\geq 1\) be a large number, and let \(\varphi (n)\) be the Euler totient function. Then
	\begin{equation} \label{el88500}
	\sum_{p\leq x} \frac{1}{p}\sum_{\gcd(n,p-1)=1} 1=a_1\li(x)+O\left(\frac{x}{(\log x)^c}\right) ,
	\end{equation}
where \(c> 1\) is an arbitrary constant, as \(x \rightarrow \infty\).  
\end{lem}

\begin{proof}[\textbf{Proof}] A routine rearrangement and an application of \hyperlink{lem7776.200}{Lemma} \ref{lem7776.200} give 
\begin{eqnarray} \label{eq7776.200h}
\sum_{p\leq x} \frac{1}{p}\sum_{\gcd(n,p-1)=1} 1&=&\sum_{p\leq x} \frac{\varphi(p-1)}{p-1}-\sum_{p\leq x} \frac{\varphi(p-1)}{p(p-1)}\\
&=&a_1 \li(x)+O \left (\frac{x}{(\log x)^c}\right ) \nonumber,
\end{eqnarray} 
where $c>1$ is an arbitrary constant. 
\end{proof}

\section{Estimate for the Error Term} \label{S6-S}\hypertarget{S6-S}
The upper bounds for exponential sums over subsets of elements in finite fields $\mathbb{F}_p$ studied in \hyperlink{S9933Q-A}{Section} \ref{S9933Q-A} are used to estimate the error term $E(x)$ in the proof of \hyperlink{thm7776.100}{Theorem} \ref{thm7776.100}. 

\begin{lem} \label{lem6.1}\hypertarget{lem6.1}  Let \(p\geq 2\) be a large prime, let \(\psi \neq 1\) be an additive character, and let \(\tau\) be a primitive root mod \(p\). Assume the element \(u\ne 0\) is not a primitive root. Then, 
	\begin{equation} \label{el89400}
	\sum_{x \leq p\leq 2x}
	\frac{1}{p}\sum_{\substack{1\leq n\leq p-1\\\gcd(n,p-1)=1}} \sum_{ 0<k\leq p-1} \psi \left((\tau ^n-u)k\right)\ll x^{1/2+\delta}\log x 
	\end{equation} 
	for all sufficiently large numbers $x\geq 1$ and an arbitrarily small number \(\delta>0\).
\end{lem}

\begin{proof}[\textbf{Proof}]  Let $\psi(z)=e^{i 2 \pi z/p}$. By hypothesis $u \ne \tau^n$ for any $n \geq 1$ such that $\gcd(n,p-1)=1$. Rearrange the triple finite sum in the form of a discrete Fourier transform
\begin{eqnarray} \label{eq89901}
E(x)&=&\sum_{x \leq p \leq 2x}\frac{1}{p} \sum_{ 0<k\leq p-1,}  \sum_{\substack{1\leq n\leq p-1\\\gcd(n,p-1)=1}} \psi ((\tau ^n-u)k) \\  
&= & \sum_{x \leq p \leq 2x}\frac{1}{p}
 \left (\sum_{ 0<k\leq p-1} e^{ \frac{-i 2 \pi uk}{p}}  \sum_{\substack{1\leq n\leq p-1\\\gcd(n,p-1)=1}} e^{ \frac{i 2 \pi k\tau ^n}{p}} \right )\nonumber.
\end{eqnarray} 
Applying \hyperlink{lem9933ERP.230V}{Lemma} \ref{lem9933ERP.230V} to the inner double sum 
and completing the upper bound of the finite sum return
\begin{eqnarray} \label{el89991}
|E(x)|	&\ll & \sum_{x \leq p \leq 2x} \frac{ 1}{p}  \cdot p^{1/2+\delta} (\log p)^2\\[.2cm]
	&\ll &  \frac{ (\log x)^2}{x^{1/2-\delta }}\sum_{x \leq p \leq 2x}  1 \nonumber \\[.3cm]
	&\ll &  x^{1/2+\delta}\log x \nonumber,
\end{eqnarray}
where the number of primes in the short interval $[x,2x]$ is $\pi(2x)-\pi(x)\leq 2 x/ \log x$. 

\end{proof}

\section{Proof of Main Result} \label{S7777-K}\hypertarget{S7777-K}
The counting function for the number of primes $p\in [x,2x]$, with a fixed primitive $u=\pm 1,v^2$, is defined by
\begin{equation} \label{el700}
\pi(u,x)=\#\{x\leq p\leq 2x: \ord_p u= p-1\}.
\end{equation}

\begin{proof}[\textbf{Proof}] (\hyperlink{thm7776.100}{Theorem} \ref{thm7776.100}): Suppose that \(u\neq \pm 1,v^2\) is not a primitive root for all primes \(p\geq x_0\), with \(x_0\geq 1\) constant. Let \(x>x_0\) be a large number, and consider the sum of the characteristic function over the short interval \([x,2x]\), that is, \\
\begin{equation} \label{el88740}
\pi_u(x)=\sum _{x \leq p\leq 2x} \Psi (u)=0.
\end{equation}\\
Replacing the characteristic function, \hyperlink{lem9955FF.200DF}{Lemma} \ref{lem9955FF.200DF}, and expanding the nonexistence equation (\ref{el88740}) yield\\
\begin{eqnarray} \label{el88750}
\pi_u(x)&=&\sum _{x \leq p\leq 2x} \Psi (u)  \\[.2cm]
&=&\sum_{x \leq p\leq 2x} \frac{1}{p}\sum_{\gcd(n,p-1)=1,} \sum_{ 0\leq k\leq p-1} \psi \left((\tau ^n-u)k\right)\nonumber\\[.2cm]
&=&\delta(u)\sum_{x \leq p\leq 2x} \frac{1}{p}\sum_{\gcd(n,p-1)=1} 1+\sum_{x \leq p\leq 2x}
\frac{1}{p}\sum_{\gcd(n,p-1)=1,} \sum_{ 0<k\leq p-1} \psi \left((\tau ^n-u)k\right)\nonumber\\[.4cm]
&=&\delta(u)M(x) + E(x)\nonumber,
\end{eqnarray} \\
where $\delta(u) \geq0$ is a constant depending on the fixed integer $u\ne 0$. \\

The main term $M(x)$ is determined by a finite sum over the trivial additive character \(\psi =1\), and the error term $E(x)$ is determined by a finite sum over the nontrivial additive characters \(\psi(t) =e^{i 2\pi t  /p}\neq 1\).\\

Applying \hyperlink{lem7776.200B}{Lemma} \ref{lem7776.200B} to the main term, and \hyperlink{lem6.1}{Lemma} \ref{lem6.1} to the error term yield\\
\begin{eqnarray} \label{el89760}
\sum _{x \leq p\leq 2x} \Psi (u)
&=&\delta(u)M(x) + E(x) \\[.2cm]
&=&\delta(u)\left (\text{li}(2x)-\text{li}(x) \right )+O\left(\frac{x}{(\log x)^B}\right)+O\left (x^{1/2+\delta}\log x  \right )\nonumber\\[.3cm]
&=& \delta(u)\frac{x}{\log x}+O\left(\frac{x}{(\log x)^B}\right) \nonumber \\[.3cm]
&>&0 \nonumber,
\end{eqnarray} 
where $B>1$. However, $\delta(u) > 0$ contradicts the hypothesis  (\ref{el88740}) for all sufficiently large numbers $x \geq x_0$. Ergo, the short interval $[x,2x]$ contains primes with the fixed primitive root $u$.   
\end{proof}

A formula for computing the density constant \(\delta(u) \geq 0\), which is the density of primes with a fixed primitive root \(u\ne \pm 1,v^2\), is derived in {\color{red}\cite[p. 220]{HC1967}}. This formula specifies the density as follows. Let $u=(st^2)^k \ne \pm 1, v^2$ with $s$ squarefree, and $k \geq 1$, and let 
\begin{equation} \label{666a}
a_k(u)=
\prod_{p \mid k }\frac{1}{p-1} \prod_{p \nmid k} \left ( 1 -\frac{1}{p(p-1)} \right) .
\end{equation}

\begin{enumerate}
\item If $s \not \equiv 1 \bmod 4$, then, 
\begin{equation} \label{666b}
\delta(u)=a_k(u). 
\end{equation}
\item If $s \equiv 1 \bmod 4$, then,  
\begin{equation} \label{666c}
\delta(u)=
\left (1-\mu(s) \prod_{\substack{p \mid  s \\ p \mid k }}\frac{1}{p-2}\prod_{\substack{p\mid s \\ p \nmid k }}\frac{1}{p^2-p-1} \right )
a_k(u). 
\end{equation}
\end{enumerate}

In the case of a decimal base $u=10$, the average density \eqref{eq7776.101d}, which is computed in \cite{WJ1961}, is
\begin{equation} \label{667}
\delta(10)= \prod_{p \geq 2} \left ( 1 -\frac{1}{p(p-1)} \right) =0.37739558136192022880547280\ldots .
\end{equation}\\

The argument can be repeated for any other squarefree base $u=2, 3, \ldots $.

Information on the determination of the constant \(\delta(u) \geq 0\), which is the density of primes with a fixed primitive root \(u\ne \pm 1,v^2\), is provided in the next section and in {\color{red}\cite[p. 220]{HC1967}}. The range of values of the density 
\begin{align}
\frac{84}{85}a_1\leq \delta(u)\leq \frac{6}{5}a_1
\end{align} is studied in \cite{PS2024}, where the constant $a_1$ is defined in \eqref{eq7776.200f}.
\section{Application To Repeated Decimals} \label{S8-R}\hypertarget{S8-R}
Let \(p\geq 2\) be a prime. The period of the repeating decimal number
\begin{equation} 
1/p=0.\overline{x_{d-1} \ldots x_1x_0},
\end{equation} 
with $x_i \in \{0,1, 2, \ldots , 9\}$, was investigated by Gauss and earlier authors centuries ago, see \cite{BM2009} for a historical account, and \cite{MM1988} for recent developments. As discussed in Articles 14-18 in \cite{GC1986}, the period, denoted by $\ord_p(10)=d \geq 1$, is a divisor of $p-1$. The problem of computing the densities for the subsets of primes for which the repeating decimals have very large periods such as $d=(p-1)/2$, and $d=p-1$, is a recent problem. This note considers the following result.\\

\begin{thm} \label{thm8.1}\hypertarget{thm8.1} There are infinitely many primes \(p\geq 7\) with maximal repeating decimal 
\begin{equation} \label{el03b}
1/p=0.\overline{x_{p-2}x_{p-3} \cdots x_1x_0},
\end{equation} 
where $ 0 \leq x_i\leq 9$. Moreover, the counting function for these primes satisfies the lower bound 
\begin{equation} \label{el03}
\pi _{10}(x)=\#\left\{ p\leq x:\ord_p(10)=p-1 \right\} =\delta(10)\frac{x}{\log x}+O\left(\frac{x}{(\log x)^2} \right),
\end{equation}
where $\delta(10)=0.3773955 \ldots$ is a constant, see {\normalfont (\ref{667})}, for all large numbers \(x\geq 1\).
\end{thm} 

The proof of \hyperlink{thm8.1}{Theorem} \ref{thm8.1} is a corollary of the more general result in \hyperlink{thm7776.100}{Theorem} \ref{thm7776.100} presented in \hyperlink{S7777-K}{Section} \ref{S7777-K}. This analysis generalizes to repeating $\ell$-adic expansions 
\begin{equation} \label{el03c}
1/p=0.\overline{x_{d-1}x_{d-2} \cdots x_1x_0},
\end{equation}
where $ 0 \leq x_i\leq \ell-1$, in any numbers system with nonsquare integer base $\ell \geq 2$.

\subsection{Proof of Theorem 8.1} \label{S9-F}\hypertarget{S9-F}
The repeating decimal fractions have the squarefree base $u=10$. In particular, the repeated fraction representation 
\begin{equation}
\frac{1}{p}=\frac{m}{10^d}+\frac{m}{10^{2d}}+ \cdots= m \sum_{n \geq 1} \frac{1}{10^{dn}}=\frac{m}{10^d-1}
\end{equation}
has the maximal period $d=p-1$ if and only if 10 has order $\ord_p(10)=p-1$ modulo $p$. This follows from the Fermat little theorem and $10^d-1 =mp$. The subset of Abel-Wieferich primes, which satisfy the congruence 
\begin{equation}\label{8825}
10^{p-1}-1\equiv 0 \bmod p^2,
\end{equation}
are also important in the calculation of the period modulo $p^2$, confer {\color{red}\cite[p.\ 333]{RP1996}}, \cite{DK2011} for other details. \\ 

The result for repeating decimal of maximal period is a simple corollary of the previous result. \\

\begin{proof}[\textbf{Proof}] \text{(\hyperlink{thm8.1}{Theorem} \ref{thm8.1})} This follows from \hyperlink{thm7776.100}{Theorem} \ref{thm7776.100} replacing the base $u=10$. That is,
\begin{equation} \label{el8964}
\pi_{10}(x)=\sum _{p\leq x} \Psi (10)
=\delta(10)\li(x) +O\left( \frac{x}{(\log x)^B } \right) ,
\end{equation}    
where $B>1$ is a constant. 
\end{proof}


\end{document}